\documentclass{amsart}

\usepackage{amsmath}
\usepackage{amsfonts}
\usepackage{amssymb,enumerate}
\usepackage{amsthm}
\usepackage[all]{xy}
\usepackage{hyperref}

\DeclareMathOperator{\hight}{ht}

\newcommand{\Spec}{\operatorname{Spec}}
\newcommand{\QSpec}{\operatorname{QSpec}}

\renewcommand{\dim}{\operatorname{dim}}

\newcommand{\Na}{\operatorname{Na}}
\newcommand{\Kr}{\operatorname{Kr}}

\newcommand{\Max}{\operatorname{Max}}
\newcommand{\QMax}{\operatorname{QMax}}

\newcommand{\fQ}{\frak{Q}}

\newtheorem{thm}{Theorem}[section]
\newtheorem{cor}[thm]{Corollary}
\newtheorem{lem}[thm]{Lemma}
\newtheorem{prop}[thm]{Proposition}
\newtheorem{defn}[thm]{Definition}
\newtheorem{exam}[thm]{Example}
\newtheorem{rem}[thm]{Remark}
\newtheorem{ques}[thm]{Question}

\begin{document}

\bibliographystyle{amsplain}

\date{}

\author{Parviz Sahandi}

\address{Department of Mathematics, University of Tabriz, Tabriz, Iran and School of Mathematics, Institute for
Research in Fundamental Sciences (IPM), Tehran Iran.}

\email{sahandi@ipm.ir, sahandi@tabrizu.ac.ir}



\title{Semistar-Krull and valuative dimension of integral domains}

\begin{abstract} Given a stable semistar operation of finite type $\star$ on
an integral domain $D$, we show that it is possible to define in a
canonical way a stable semistar operation of finite type $\star[X]$
on the polynomial ring $D[X]$, such that, if $n:=\star$-$\dim(D)$,
then $n+1\leq \star[X]\text{-}\dim(D[X])\leq 2n+1$. We also
establish that if $D$ is a $\star$-Noetherian domain or is a
Pr\"{u}fer $\star$-multiplication domain, then
$\star[X]\text{-}\dim(D[X])=\star\text{-}\dim(D)+1$. Moreover we
define the semistar valuative dimension of the domain $D$, denoted
by $\star$-$\dim_v(D)$, to be the maximal rank of the
$\star$-valuation overrings of $D$. We show that
$\star$-$\dim_v(D)=n$ if and only if
$\star[X_1,\cdots,X_n]$-$\dim_v(D[X_1,\cdots,X_n])=2n$, and that if
$\star$-$\dim_v(D)<\infty$ then
$\star[X]$-$\dim_v(D[X])=\star$-$\dim_v(D)+1$. In general
$\star$-$\dim(D)\leq\star$-$\dim_v(D)$ and equality holds if $D$ is
a $\star$-Noetherian domain or is a Pr\"{u}fer
$\star$-multiplication domain. We define the $\star$-Jaffard domains
as domains $D$ such that $\star$-$\dim(D)<\infty$ and
$\star$-$\dim(D)=\star$-$\dim_v(D)$. As an application,
$\star$-quasi-Pr\"{u}fer domains are characterized as domains $D$
such that each $(\star,\star')$-linked overring $T$ of $D$, is a
$\star'$-Jaffard domain, where $\star'$ is a stable semistar
operation of finite type on $T$. As a consequence of this result we
obtain that a Krull domain $D$, must be a $w_D$-Jaffard domain.
\end{abstract}

\maketitle

\section{Introduction}

\noindent Throughout this paper, $D$ denotes a (commutative
integral) domain with identity and $K$ denotes the quotient field of
$D$. Let $X$ be an algebraically independent indeterminate over $D$.
Seidenberg proved in \cite[Theorem 2]{S1}, that if $D$ has finite
Krull dimension, then
$$\dim(D)+1\leq\dim(D[X])\leq2(\dim(D))+1.$$
Moreover, Krull \cite{Kr} has shown that if $D$ is any
finite-dimensional Noetherian ring, then $\dim(D[X])=1+\dim(D)$ (cf.
also \cite[Theorem 9]{S1}). Seidenberg subsequently proved the same
equality in case $D$ is any finite-dimensional Pr\"{u}fer domain. To
unify and extend such results on Krull-dimension, Jaffard
\cite{Jaf1} introduced and studied the \emph{valuative dimension}
denoted by $\dim_v(D)$, for a domain $D$. This is the maximum of the
ranks of the valuation overrings of $D$. Jaffard proved in
\cite[Chapitre IV]{Jaf1} (see also Arnold \cite{A}), that if $D$ has
finite valuative dimension, then $\dim_v(D[X])=1+\dim_v(D)$ and that
if $D$ is a Noetherian or a Pr\"{u}fer domain, then
$\dim(D)=\dim_v(D)$. Also he showed that $\dim_v(D)=n$ if and only
if $\dim(D[X_1,\cdots,X_n])=2n$, where $X_1,\cdots,X_n$ are
indeterminates over $D$. In \cite{ABDFK} the authors introduced the
notion of Jaffard domains, as integral domains $D$ such that
$\dim(D)=\dim_v(D)$. The class of Jaffard domains contains most of
the well-known classes of finite dimensional rings involved in
dimension theory of commutative rings, such as Noetherian domains,
Pr\"{u}fer domains, universally catenarian domains \cite{BDF}, and
stably strong S-domains \cite{MM, Kab}. A good and available
reference for the dimension theory of commutative rings is Gilmer
\cite[Section 30]{G}.

For several decades, star operations, as described in \cite[Section
32]{G}, have proven to be an essential tool in \emph{multiplicative
ideal theory}, for studying various classes of domains. In
\cite{OM}, Okabe and Matsuda introduced the concept of a semistar
operation to extend the notion of a star operation. Since then,
semistar operations have been extensively studied and, because of a
greater flexibility than star operations, have permitted a finer
study and new classifications of special classes of integral
domains. For instance, semistar-theoretic analogues of the classical
notions of Krull dimension, Noetherian and Pr\"{u}fer domains have
been introduced: see \cite{EFP} and \cite{FJS} for the basics on
$\star$-Krull dimension, $\star$-Noetherian domains and Pr\"{u}fer
$\star$-multiplication domains (for short P$\star$MD), respectively.

Now it is natural to ask:

\begin{ques}\label{1111} Given a semistar operation of finite type $\star$ on $D$, is it
possible to define in a canonical way a semistar operation of finite
type $\star[X]$ on $D[X]$, such that $\star\text{-}
\dim(D)+1\leq\star[X]\text{-}\dim(D[X])\leq2(\star\text{-}\dim(D))+1$,
and that if $D$ is a $\star$-Noetherian domain or a P$\star$MD, then
$\star[X]\text{-}\dim(D[X])=\star\text{-}\dim(D)+1$? \end{ques}

In this paper, we answer this question, in case that $\star$ is a
stable semistar operation of finite type on $D$. More precisely, in
Section 2, using the technique introduced by Chang and Fontana in
\cite{CF1}, we define in a canonical way a semistar operation stable
and of finite type $\star[X]$ on $D[X]$: see Theorem \ref{main}. In
Section 3 we show among other things that this question has an
affirmative answer: see Theorems \ref{dim}, \ref{nodim}, and
\ref{prudim}.

Let $\star$ be a semistar operation on the integral domain $D$ and
let $\widetilde{\star}$ be the stable semistar operation of finite
type canonically associated to $\star$ (the definitions are recalled
later in this section). We define in Section 4, what it means the
semistar valuative dimension of $D$, denoted by $\star$-$\dim_v(D)$.
It extends the ``classical'' valuative dimension of P. Jaffard
\cite{Jaf1}, denoted by $\dim_v(D)$ to the setting of semistar
operations. We show that the semistar valuative dimension of $D$ has
various nice properties, like the classical valuative dimension. For
example we show that if $\widetilde{\star}$-$\dim_v(D)<\infty$ then
$\star[X]\text{-}\dim_v(D[X])=\widetilde{\star}\text{-}\dim_v(D)+1$:
see Theorem \ref{VV}. Also we established that
$\widetilde{\star}\text{-}\dim(D)\leq\widetilde{\star}\text{-}\dim_v(D)$,
and equality holds if $D$ is a $\widetilde{\star}$-Noetherian domain
or a P$\star$MD: see Corollaries \ref{novdim} and \ref{prdim}. In
relation with the $\star$-Nagata ring $\Na(D,\star)$, it is shown
that $\widetilde{\star}\text{-}\dim_v(D)=\dim_v(\Na(D,\star))$: see
Theorem \ref{vdim}. If $\widetilde{\star}\text{-}\dim(D)<\infty$ and
$\widetilde{\star}\text{-}\dim(D)=\widetilde{\star}\text{-}\dim_v(D)$,
we say that, $D$ is a \textit{$\widetilde{\star}$-Jaffard domain}.
We establish that $D$ is a $\widetilde{\star}$-quasi-Pr\"{u}fer
domain if and only if each $(\star,\star')$-linked overring $T$ of
$D$ is a $\widetilde{\star'}$-Jaffard domain, where $\star'$ is a
semistar operation on $T$: see Theorem \ref{qJaf}. As a consequence
of this result we obtain that a Krull domain $D$, must be a
$w_D$-Jaffard domain.

To facilitate the reading of the introduction and of the paper, we
first review some basic facts on semistar operations. Let
$\overline{\mathcal{F}}(D)$ denote the set of all nonzero
$D$-submodules of $K$. Let $\mathcal{F}(D)$ be the set of all
nonzero \emph{fractional} ideals of $D$; i.e., $E\in\mathcal{F}(D)$
if $E\in\overline{\mathcal{F}}(D)$ and there exists a nonzero
element $r\in D$ with $rE\subseteq D$. Let $f(D)$ be the set of all
nonzero finitely generated fractional ideals of $D$. Obviously,
$f(D)\subseteq\mathcal{F}(D)\subseteq\overline{\mathcal{F}}(D)$. As
in \cite{OM}, a {\it semistar operation on} $D$ is a map
$\star:\overline{\mathcal{F}}(D)\rightarrow\overline{\mathcal{F}}(D)$,
$E\mapsto E^{\star}$, such that, for all $x\in K$, $x\neq 0$, and
for all $E, F\in\overline{\mathcal{F}}(D)$, the following three
properties hold:
\begin{itemize}
\item [$\star_1$]:  $(xE)^{\star}=xE^{\star}$;
\item [$\star_2$]:  $E\subseteq F$ implies that $E^{\star}\subseteq
F^{\star}$;
\item [$\star_3$]:  $E\subseteq E^{\star}$ and
$E^{\star\star}:=(E^{\star})^{\star}=E^{\star}$.
\end{itemize}
Recall from \cite[Proposition 5]{OM} that if $\star$ is a semistar
operation on $D$, then, for all $E, F\in\overline{\mathcal{F}}(D)$,
the following basic formulas follow easily from the above axioms:
\begin{itemize}
\item [(1)] $(EF)^{\star}=(E^{\star}F)^{\star}=(EF^{\star})^{\star}=(E^{\star}F^{\star})^{\star}$;
\item [(2)] $(E+F)^{\star}=(E^{\star}+F)^{\star}=(E+F^{\star})^{\star}=(E^{\star}+F^{\star})^{\star}$;
\item [(3)] $(E:F)^{\star}\subseteq(E^{\star}:F^{\star})=(E^{\star}:F)=(E^{\star}:F)^{\star}$, if
$(E:F)\neq(0)$;
\item [(4)] $(E\cap F)^{\star}\subseteq E^{\star}\cap F^{\star}=(E^{\star}\cap F^{\star})^{\star}$ if $(E\cap F)\neq(0)$.
\end{itemize}

It is convenient to say that a \emph{(semi)star operation on} $D$ is
a semistar operation which, when restricted to $\mathcal{F}(D)$, is
a star operation (in the sense of \cite[Section 32]{G}). It is easy
to see that a semistar operation $\star$ on $D$ is a (semi)star
operation on $D$ if and only if $D^{\star}=D$.

Let $\star$ be a semistar operation on the domain $D$. For every
$E\in\overline{\mathcal{F}}(D)$, put $E^{\star_f}:=\cup F^{\star}$,
where the union is taken over all finitely generated $F\in f(D)$
with $F\subseteq E$. It is easy to see that $\star_f$ is a semistar
operation on $D$, and ${\star_f}$ is called \emph{the semistar
operation of finite type associated to} $\star$. Note that
$(\star_f)_f=\star_f$. A semistar operation $\star$ is said to be of
\emph{finite type} if $\star=\star_f$; in particular ${\star_f}$ is
of finite type. We say that a nonzero ideal $I$ of $D$ is a
\emph{quasi-$\star$-ideal} of $D$, if $I^{\star}\cap D=I$; a
\emph{quasi-$\star$-prime} (ideal of $D$), if $I$ is a prime
quasi-$\star$-ideal of $D$; and a \emph{quasi-$\star$-maximal}
(ideal of $D$), if $I$ is maximal in the set of all proper
quasi-$\star$-ideals of $D$. Each quasi-$\star$-maximal ideal is a
prime ideal. It was shown in \cite[Lemma 4.20]{FH} that if
$D^{\star} \neq K$, then each proper quasi-$\star_f$-ideal of $D$ is
contained in a quasi-$\star_f$-maximal ideal of $D$. We denote by
$\QMax^{\star}(D)$ (resp., $\QSpec^{\star}(D)$) the set of all
quasi-$\star$-maximal ideals (resp., quasi-$\star$-prime ideals) of
$D$. When $\star$ is a (semi)star operation, it is easy to see that
the notion of quasi-$\star$-ideal is equivalent to the classical
notion of $\star$-ideal (i.e., a nonzero ideal $I$ of $D$ such that
$I^{\star}=I$).

If $\star_1$ and $\star_2$ are semistar operations on $D$, one says
that $\star_1\leq\star_2$ if $E^{\star_1}\subseteq E^{\star_2}$ for
each $E\in\overline{\mathcal{F}}(D)$ (cf. \cite[page 6]{OM}). This
is equivalent to saying that
$(E^{\star_1})^{\star_2}=E^{\star_2}=(E^{\star_2})^{\star_1}$ for
each $E\in\overline{\mathcal{F}}(D)$ (cf. \cite[Lemma 16]{OM}).
Obviously, for each semistar operation $\star$ defined on $D$, we
have $\star_f\leq\star$. Let $d_D$ (or, simply, $d$) denote the
identity (semi)star operation on $D$. Clearly, $d_D\leq\star$ for
all semistar operations $\star$ on $D$.

If $\Delta$ is a set of prime ideals of a domain $D$, then there is
an associated semistar operation on $D$, denoted by
$\star_{\Delta}$, defined as follows:
$$E^{\star_{\Delta}}:=\cap\{ED_P|P\in\Delta\}\text{, for each }E\in\overline{\mathcal{F}}(D).$$
If $\Delta=\emptyset$, let $E^{\star_{\Delta}}:=K$ for each
$E\in\overline{F}(D)$. Note that $E^{\star_{\Delta}}D_P=ED_P$ for
each $E\in\overline{\mathcal{F}}(D)$ and $P\in\Delta$ by \cite[Lemma
4.1 (2)]{FH}. One calls $\star_{\Delta}$ the \emph{spectral semistar
operation associated to} $\Delta$. A semistar operation $\star$ on a
domain $D$ is called a \emph{spectral semistar operation} if there
exists a subset $\Delta$ of the prime spectrum of $D$, Spec$(D)$,
such that $\star=\star_{\Delta}$. When $\Delta:=\QMax^{\star_f}(D)$,
we set $\widetilde{\star}:=\star_{\Delta}$; i.e.,
$$E^{\widetilde{\star}}:= \cap\{ED_P|P\in\QMax^{\star_f}(D)\}\text{, for each
}E\in\overline{\mathcal{F}}(D).$$

It has become standard to say that a semistar operation $\star$ is
{\it stable} if $(E\cap F)^{\star}=E^{\star}\cap F^{\star}$ for all
$E$, $F\in \overline{\mathcal{F}}(D)$. (``Stable" has replaced the
earlier usage, ``quotient", in \cite[Definition 21]{OM}.) All
spectral semistar operations are stable \cite[Lemma 4.1(3)]{FH}. In
particular, for any semistar operation $\star$, we have that
$\widetilde{\star}$ is a stable semistar operation of finite type
\cite[Corollary 3.9]{FH}.

Let $D$ be a domain, $\star$ a semistar operation on $D$, $T$ an
overring of $D$, and $\iota:D\hookrightarrow T$ the corresponding
inclusion map. In a canonical way, one can define an associated
semistar operation $\star_{\iota}$ on $T$, by $E\mapsto
E^{\star_{\iota}}:=E^{\star}$, for each
$E\in\overline{\mathcal{F}}(T)(\subseteq\overline{\mathcal{F}}(D))$.

The most widely studied (semi)star operations on $D$ have been the
identity $d_D$ and $v_D$, $t_D:=(v_D)_f$, and $w_D:=\widetilde{v_D}$
operations, where $E^{v_D}:=(E^{-1})^{-1}$, with
$E^{-1}:=(D:E):=\{x\in K|xE\subseteq D\}$.

Let $D$ be a domain with quotient field $K$, and let $X$ be an
indeterminate over $K$. For each $f\in K[X]$, we let $c_D(f)$ denote
the content of the polynomial $f$, i.e., the (fractional) ideal of
$D$ generated by the coefficients of $f$. Let $\star$ be a semistar
operation on $D$. If $N_{\star}:=\{g\in D[X]|g\neq0\text{ and
}c_D(g)^{\star}=D^{\star}\}$, then $N_{\star}=
D[X]\backslash\bigcup\{P[X]|P\in\QMax^{\star_f}(D)\}$ is a saturated
multiplicative subset of $D[X]$. The ring of fractions
$$\Na(D,\star):=D[X]_{N_{\star}}$$ is called the $\star$-{\it Nagata domain (of $D$ with respect to the
semistar operation} $\star$). When $\star=d$, the identity
(semi)star operation on $D$, then $\Na(D,d)$ coincides with the
classical Nagata domain $D(X)$ (as in, for instance \cite[page
18]{Na}, \cite[Section 33]{G} and \cite{FL}).

\section{Semistar operations on polynomial rings}

In \cite{CF1}, Chang and Fontana introduced a new technique for
defining new semistar operations on integral domains. Let $D$ be an
integral domain with quotient field $K$, and let $X$ be an
indeterminate over $K$. For a given multiplicative subset
$\mathcal{S}$ of $D[X]$, set
$$E^{\circlearrowleft_{\mathcal{S}}}:=E[X]_{\mathcal{S}}\cap
K,\text{ for all }E\in \overline{\mathcal{F}}(D).$$ Then it is
proved in \cite[Theorem 2.1]{CF1} among other things that, the
mapping
$\circlearrowleft_{\mathcal{S}}:\overline{\mathcal{F}}(D)\to\overline{\mathcal{F}}(D)$,
$E\mapsto E^{\circlearrowleft_{\mathcal{S}}}$ is a stable semistar
operation of finite type on $D[X]$, i.e.,
$\widetilde{\circlearrowleft_{\mathcal{S}}}=\circlearrowleft_{\mathcal{S}}$,
and $\QMax^{\circlearrowleft_{\mathcal{S}}}(D)=$ the set of maximal
elements of
$\Delta(\mathcal{S}):=\{P\in\Spec(D)|P[X]\cap\mathcal{S}=\emptyset\}$.

Let $D$ be an integral domain, and $\star$ a semistar operation on
$D$. Using the technique discussed in the first paragraph, Chang and
Fontana defined canonically a semistar operation denoted by
$[\star]$ on the polynomial ring $D[X]$. More precisely suppose that
$X$, $Y$ are two indeterminates over $D$, and set $D_1:=D[X]$,
$K_1:=K(X)$. Take the following subset of $\Spec(D_1)$:
$$\Delta_1^{\star}:=\{Q_1\in\Spec(D_1)|\text{ }Q_1\cap D=(0)\text{ or }
Q_1=(Q_1\cap D)[X]\text{ and }(Q_1\cap D)^{\star_f}\subsetneq
D^{\star}\}.$$ Set
$\mathcal{S}_1^{\star}:=\mathcal{S}(\Delta_1^{\star}):=D_1[Y]\backslash(\bigcup\{Q_1[Y]
|Q_1\in\Delta_1^{\star}\})$ and
$[\star]:=\circlearrowleft_{\mathcal{S}_1^{\star}}$, that is:
$$E^{[\star]}:=E[Y]_{\mathcal{S}_1^{\star}}\cap
K_1, \text{   for all }E\in \overline{\mathcal{F}}(D_1).$$ They
proved answering their question \cite[Question]{CF}, that $D$ is a
$\widetilde{\star}$-quasi-Pr\"{u}fer domain if and only if each
upper to zero, is a quasi-$[\star]$-maximal ideal of $D[X]$. Recall
that $D$ is said to be a \textit{$\star$-quasi-Pr\"{u}fer domain},
in case, if $Q$ is a prime ideal in $D[X]$, and $Q\subseteq P[X]$,
for some $P\in\QSpec^{\star}(D)$, then $Q=(Q\cap D)[X]$. This notion
is the semistar analogue of the classical notion of the
quasi-Pr\"{u}fer domains \cite[Section 6.5]{FHP} (that is among
other equivalent conditions, the domain $D$ is said to be a
\emph{quasi-Pr\"{u}fer domain} if it has Pr\"{u}ferian integral
closure).

Now by the same technique, we define canonically a semistar
operation denoted by $\star[X]$ on the polynomial ring $D[X]$, which
has desired semistar (Krull) dimension theoretic properties.

\begin{thm}\label{main} Let $D$ be an integral domain with quotient field $K$, let $X$, $Y$ be two indeterminates
over $D$ and let $\star$ be a semistar operation on D. Set
$D_1:=D[X]$, $K_1:=K(X)$ and take the following subset of
$\Spec(D_1)$:
$$\Theta_1^{\star}:=\{Q_1\in\Spec(D_1)|\text{ }Q_1\cap D=(0)\text{ or }(Q_1\cap D)^{\star_f}\subsetneq D^{\star}\}.$$
Set
$\mathfrak{S}_1^{\star}:=\mathcal{S}(\Theta_1^{\star}):=D_1[Y]\backslash(\bigcup\{Q_1[Y]
|Q_1\in\Theta_1^{\star}\})$ and:
$$E^{\circlearrowleft_{\mathfrak{S}_1^{\star}}}:=E[Y]_{\mathfrak{S}_1^{\star}}\cap
K_1, \text{   for all }E\in \overline{\mathcal{F}}(D_1).$$

\begin{itemize}
\item[(a)] The mapping $\star[X]:=\circlearrowleft_{\mathfrak{S}_1^{\star}}:
\overline{\mathcal{F}}(D_1)\to\overline{\mathcal{F}}(D_1)$,
$E\mapsto E^{\circlearrowleft_{\mathfrak{S}_1^{\star}}}$ is a stable
semistar operation of finite type on $D[X]$, i.e.,
$\widetilde{\star[X]}=\star[X]$.

\item[(b)] $\widetilde{\star}[X]=\star_f[X]=\star[X]$.

\item[(c)] $\star[X]\leq[\star]$. In particular, if $\star$ is a (semi)star
operation on $D$, then $\star[X]$ is a (semi)star operation on
$D[X]$.

\item[(d)] $d_D[X]=d_{D[X]}$.
\end{itemize}
\end{thm}

\begin{proof} Note that, if $Q_1\in\Spec(D[X])$ is not an upper to zero
and $(Q_1\cap D)^{\star_f}\subsetneq D^{\star}$, then the prime
ideal $Q_1\cap D$ is contained in a quasi-$\star_f$-maximal ideal of
$D$. Moreover if $Q_1\cap D=(0)$ and $c_D(Q_1)^{\star_f}\subsetneq
D^{\star}$ then $c_D(Q_1)^{\star_f}$ is contained in a
quasi-$\star_f$-prime ideal $P$ of $D$ and hence $Q_1\subseteq P[X]$
with $P^{\star_f}\subsetneq D^{\star}$. Set
$\circleddash_1^{\star}:=\{Q_1\in\Spec(D_1)|\text{ }Q_1\cap
D=(0)\text{ and }c_D(Q_1)^{\star_f}=D^{\star}\text{ or }Q_1\cap
D\in\QMax^{\star_f}(D)\}.$ Now we show that:
$$\mathfrak{S}_1^{\star}:=D_1[Y]\backslash(\bigcup\{Q_1[Y]
|Q_1\in\Theta_1^{\star}\})=D_1[Y]\backslash(\bigcup\{Q_1[Y]
|Q_1\in\circleddash_1^{\star}\})=\mathcal{S}(\circleddash_1^{\star}).$$
Since $\circleddash_1^{\star}\subseteq\Theta_1^{\star}$ one has
$\mathfrak{S}_1^{\star}\subseteq\mathcal{S}(\circleddash_1^{\star})$.
For the other inclusion suppose that
$f\in\mathcal{S}(\circleddash_1^{\star})$. So that $f\notin Q_1[Y]$
for each $Q_1\in\circleddash_1^{\star}$. We want to show that
$f\in\mathfrak{S}_1^{\star}$. Suppose the contrary, hence $f\in
Q_1[Y]$ for some $Q_1\in\Theta_1^{\star}$. Therefore there are two
cases to consider:

1) If $Q_1\cap D=(0)$, we have $c_D(Q_1)^{\star_f}\neq D^{\star}$ as
$Q_1\notin\circleddash_1^{\star}$. Thus $Q_1\subseteq P[X]$ for some
quasi-$\star_f$-prime ideal $P$ of $D$. Choose
$M\in\QMax^{\star_f}(D)$ such that $P\subseteq M$. So that
$Q_1\subseteq M[X]$ and hence $f\in M[X][Y]$ while
$M[X]\in\circleddash_1^{\star}$, which is a contradiction.

2) If $(Q_1\cap D)^{\star_f}\subsetneq D^{\star}$, then $Q_1\cap
D\subseteq M$ for some quasi-$\star_f$-maximal ideal of $D$. We have
$Q_1\cap D\neq M$, since otherwise $Q_1\in\circleddash_1^{\star}$
and $f\in Q_1[Y]$ which is a contradiction. Note that
$(Q_1+M[X])\cap D=(M[X]+(X))\cap D=M$ and $Q_1\subseteq
Q_1+M[X]\subseteq M[X]+(X)$. Therefore $f\in (M[X]+(X))[Y]$ while
$M[X]+(X)\in\circleddash_1^{\star}$ which is again a contradiction.

So that we have $f\notin Q_1[Y]$ for each $Q_1\in\Theta_1^{\star}$.
Thus $f\in\mathfrak{S}_1^{\star}$, that is
$\mathfrak{S}_1^{\star}=\mathcal{S}(\circleddash_1^{\star})$.

$(a)$ It follows from \cite[Theorem 2.1 (a) and (b)]{CF1}, that
$\star[X]$ is a stable semistar operation of finite type on $D[X]$.

$(b)$ Since $\QMax^{\star_f}(D)=\QMax^{\widetilde{\star}}(D)$, the
conclusion follows easily from the fact that
$\mathfrak{S}_1^{\widetilde{\star}}=\mathfrak{S}_1^{\star_f}=\mathfrak{S}_1^{\star}$.

$(c)$ It is easily seen that
$\mathfrak{S}_1^{\star}\subseteq\mathcal{S}_1^{\star}$. Then
$$E^{\star[X]}=E[Y]_{\mathfrak{S}_1^{\star}}\cap
K_1\subseteq(E[Y]_{\mathfrak{S}_1^{\star}})_{\mathcal{S}_1^{\star}}\cap
K_1=E[Y]_{\mathcal{S}_1^{\star}}\cap K_1=E^{[\star]}.$$ This means
that $\star[X]\leq[\star]$ by definition. Now if $\star$ is a
(semi)star operation on $D$, then $[\star]$ is a (semi)star
operation on $D[X]$ by \cite[Theorem 2.3 (a)]{CF1}. So that
$D_1\subseteq D_1^{\star[X]}\subseteq D_1^{[\star]}=D_1$, that is
$D_1^{\star[X]}=D_1$. Hence $\star[X]$ is a (semi)star operation on
$D[X]$.

$(d)$ Note that we have:
\begin{align*} \mathfrak{S}_1^{d_D}= &
D_1[Y]\backslash(\bigcup\{Q_1[Y]
|Q_1\in\Theta_1^{d_D}\}) \\[1ex]
                      = & D_1[Y]\backslash(\bigcup\{Q_1[Y] |Q_1\in\Spec(D_1)\text{ and }Q_1\cap D\neq D\})\\[1ex]
                      = & D_1[Y]\backslash(\bigcup\{Q_1[Y]
                      |Q_1\in\Max(D_1)\}).
\end{align*}
So for an element $E\in\overline{\mathcal{F}}(D_1)$ we have:
$$E=E^{d_{D[X]}}\subseteq E^{d_D[X]}=E[Y]_{\mathfrak{S}_1^{d_D}}\cap
K_1=ED_1(Y)\cap K_1=E.$$The last equality follows from
\cite[Proposition 3.4 (3)]{FL}. Thus $E^{d_{D[X]}}=E^{d_D[X]}$, that
is $d_D[X]=d_{D[X]}$.
\end{proof}

A different approach to the semistar operations on polynomial rings
is possible by using the notion of localizing system. Recall that a
\emph{localizing system of ideals} $\mathcal{F}$ of $D$ is a set of
(integral) ideals of $D$ verifying the following conditions $(a)$ if
$I\in \mathcal{F}$ and if $I\subseteq J$ , then $J\in \mathcal{F}$;
$(b)$ if $I\in \mathcal{F}$ and if $J$ is an ideal of $D$ such that
$(J:_DiD)\in\mathcal{F}$, for each $i\in I$, then $J\in\mathcal{F}$.
The relation between stable semistar operations and localizing
systems has been deeply investigated by M. Fontana and J. Huckaba in
\cite{FH} and by F. Halter-Koch in the context of module systems
\cite{FlK}. If $\star$ is a semistar operation on $D$, then
$\mathcal{F}^{\star}:=\{I$ ideal of $D|I^{\star}=D^{\star}\}$ is a
localizing system (called the \emph{localizing system associated to
$\star$}) of $D$. And if $\mathcal{F}$ is a localizing system of
$D$, then the map $E\mapsto
E^{\star_{\mathcal{F}}}:=\bigcup\{(E:J)|J\in\mathcal{F}\}$, for each
$E\in\overline{\mathcal{F}}(D)$, is a stable semistar operation on
$D$. It is proved in \cite[Proposition 3.1]{P} that if $\mathcal{F}$
is a localizing system of $D$, then $\mathcal{F}[X]:=\{A$ ideal of
$D[X]|A\cap D\in\mathcal{F}\}$ is a localizing system of the
polynomial ring $D[X]$. Now let $\star$ be a stable semistar
operation on $D$ and let $\mathcal{F}^{\star}$ be the localizing
system of $D$ associated to $\star$. Consider the localizing system
$\mathcal{F}^{\star}[X]$ of $D[X]$. Then G. Picozza \cite[Page
426]{P} introduced a semistar operation denoted by $\star'$ on the
polynomial ring $D[X]$ as $\star_{\mathcal{F}^{\star}[X]}$. He used
the semistar operation $\star'$ to provide the semistar version of
the Hilbert basis Theorem \cite[Theorem 3.3]{P}.

\begin{prop}\label{Picozza} If $\star$ is a stable semistar operation of finite type on $D$,
that is if, $\star=\widetilde{\star}$ then $\star'=\star[X]$.
\end{prop}

\begin{proof} Adapt the notation in the paragraph before the proposition.
Recall from \cite[Corollary 2.2]{CF1} that if $\mathcal{F}$ is a
localizing system of $D$, $Y$ is an indeterminate over $D$, and
$\mathcal{S}(\mathcal{F}):=D[Y]\backslash\bigcup\{Q[Y]|Q\in\Spec(D)$
and $Q\notin\mathcal{F}\}$ which is a saturated multiplicatively
closed subset of $D[Y]$, then
$\star_{\mathcal{F}}=\circlearrowleft_{\mathcal{S}(\mathcal{F})}$.
Now let $\mathcal{F}:=\mathcal{F}^{\star}[X]=\{A\text{ ideal of
}D_1|(A\cap D)^{\star}=D^{\star}\}$. Then
\begin{align*}
\mathcal{S}(\mathcal{F})= &
D_1[Y]\backslash\bigcup\{Q_1[Y]|Q_1\in\Spec(D_1)
\text{ and }Q_1\notin\mathcal{F}\} \\[1ex]
              = & D_1[Y]\backslash\bigcup\{Q_1[Y]|Q_1\in\Spec(D_1)
\text{ s.t. }Q_1\cap D=(0)\text{ or }(Q_1\cap D)^{\star}\subsetneq D^{\star}\}\\[1ex]
              = & D_1[Y]\backslash\bigcup\{Q_1[Y]|Q_1\in\Theta^{\star}_1\} \\[1ex]
              = & \mathfrak{S}^{\star}_1.
\end{align*}
Consequently
$\star'=\star_{\mathcal{F}}=\circlearrowleft_{\mathcal{S}(\mathcal{F})}
=\circlearrowleft_{\mathfrak{S}^{\star}_1}=\star[X]$
which ends the proof.
\end{proof}

Note that the semistar operation $[\star]$ has a main difference
with $\star[X]$ and $\star'$. Indeed let $\star=d_D$. Then one has
$d_D'=d_D[X]=d_{D[X]}$ by Theorem \ref{main} $(d)$ and Proposition
\ref{Picozza}. But $[d_D]\neq d_{D[X]}$. Since if $[d_D]=d_{D[X]}$,
then \cite[Corollary 2.5 (1)]{CF} implies that if $D$ is a
Pr\"{u}fer domain then $D[X]$ should be a Pr\"{u}fer domain which is
absurd.

\begin{rem}\label{q} Note that the set of quasi-$\star[X]$-prime ideals of $D[X]$, coincides
with the set $\Theta_1^{\star}\backslash\{0\}$. Indeed let $Q$ be an
element of $\Theta_1^{\star}\backslash\{0\}$. Then we have $Q[Y]\cap
\mathfrak{S}_1^{\star}=\emptyset$. Hence
\begin{align*}
Q^{\star[X]}\cap D[X]= & (Q[Y]_{\mathfrak{S}_1^{\star}}\cap K(X))\cap D[X] \\[1ex]
                      = & (Q[Y]_{\mathfrak{S}_1^{\star}}\cap D[X,Y])\cap D[X]\\[1ex]
                      = & Q[Y]\cap D[X]=Q.
\end{align*}
Therefore $Q$ is a quasi-$\star[X]$-prime ideal of $D[X]$; i.e.,
$\Theta_1^{\star}\backslash\{0\}\subseteq\QSpec^{\star[X]}(D[X])$.
Since the other inclusion is trivial, we obtain that
$\QSpec^{\star[X]}(D[X])=\Theta_1^{\star}\backslash\{0\}$.
\end{rem}

In the rest of the paper for every semistar operation $\star$ on an
integral domain $D$, we let $\star[X]$, to be the stable semistar
operation of finite type on $D[X]$ canonically associated to $\star$
as in Theorem \ref{main}(a).

Let $\star$ be a semistar operation on a domain $D$. As in
\cite{FJS} and \cite{EF} (cf. also \cite{HMM} for the case of a star
operation), $D$ is called a \emph{Pr\"{u}fer $\star$-multiplication
domain} (for short, a P$\star$MD) if each finitely generated ideal
of $D$ is $\star_f$-invertible; i.e., if
$(II^{-1})^{\star_f}=D^{\star}$ for all $I\in f(D)$. When $\star=v$,
we recover the classical notion of P$v$MD; when $\star=d_D$, the
identity (semi)star operation, we recover the notion of Pr\"{u}fer
domain.

\begin{rem} Let $\star$ be a semistar operation on an integral domain $D$. Suppose that $D[X]$ is a P$\star[X]$MD
(resp. a $\star[X]$-quasi-Pr\"{u}fer domain). Since
$\star[X]\leq[\star]$ by Theorem \ref{main}(c), we obtain that
$D[X]$ is a P$[\star]$MD by \cite{FJS}(resp. a
$[\star]$-quasi-Pr\"{u}fer domain by \cite[Corollary 2.4]{CF}). So
that $D$ is a P$\star$MD by \cite[Corollary 2.5 (1)]{CF1} (resp. a
$\widetilde{\star}$-quasi-Pr\"{u}fer domain by \cite[Corollary
2.4]{CF1}).
\end{rem}

In \cite[Section 3]{EFP}, El Baghdadi, Fontana and Picozza defined
and studied the \textit{semistar Noetherian domains}, i.e., domains
having the ascending chain condition on quasi-semistar-ideals.

\begin{rem} (Cf. \cite[Theorem 3.6]{P}) Let $\star$ be a semistar operation on an integral domain $D$. Then $D$ is a
$\widetilde{\star}$-Noetherian domain if and only if $D[X]$ is a
$\star[X]$-Noetherian domain. In fact if $D[X]$ is a
$\star[X]$-Noetherian domain, since $\star[X]\leq[\star]$ by Theorem
\ref{main}(c), we obtain that $D[X]$ is a $[\star]$-Noetherian
domain. So that $D$ is $\widetilde{\star}$-Noetherian by
\cite[Corollary 2.5 (2)]{CF1}. For the other implication use Remark
\ref{Picozza} together with \cite[Theorem 3.3]{P}.
\end{rem}

\section{Semistar-Krull dimension}

Let $\star$ be a semistar operation on an integral domain $D$. In
this section we make use of the semistar operation $\star[X]$ on
$D[X]$, canonically associated to the given semistar operation
$\star$ on $D$, to provide an answer to the question raised in the
introduction. First we recall some definitions and properties of
$\star$-dimension. For each quasi-$\star$-prime $P$ of $D$, the
\textit{$\star$-height} of $P$ (for short, $\star$-$\hight(P)$) is
defined to be the supremum of the lengths of the chains of
quasi-$\star$-prime ideals of $D$, between prime ideal $(0)$
(included) and $P$. Obviously, if $\star=d_D$ is the identity
(semi)star operation on $D$, then $\star$-$\hight(P)=\hight(P)$, for
each prime ideal $P$ of $D$. If the set of quasi-$\star$-prime of
$D$ is not empty, the $\star$-dimension of $D$ is defined as
follows:
$$\star\text{-}\dim(D):=\sup\{\star\text{-}\hight(P)|P\text{ is a
quasi-}\star\text{-prime of }D\}.$$ If the set of
quasi-$\star$-primes of $D$ is empty, then pose
$\star\text{-}\dim(D):=0$. Thus, if $\star =d_D$, then
$\star\text{-}\dim(D)= \dim(D)$, the usual (Krull) dimension of $D$.

Note that, the notions of $t$-dimension and of $w$-dimension have
received a considerable interest by several authors (cf. for
instance, \cite{F1, F2, H}).

It is known (see \cite[Lemma 2.11]{EFP}) that
\begin{align*}
\widetilde{\star}\text{-}\dim(D)= & \sup\{\hight(P) \mid P\text{ is
a quasi-}\widetilde{\star}\text{-prime ideal of } D\} \\[1ex]
                                = & \sup\{\hight(P) \mid P\text{ is a
quasi-}\widetilde{\star}\text{-maximal ideal of } D\}.
\end{align*}

We answer to the Question \ref{1111}, in the results \ref{dim},
\ref{nodim} and \ref{prdim}. The following result is the semistar
version of the classical theorem of Seidenberg \cite[Theorem 2]{S1}.

\begin{thm}\label{dim} Let $\star$ be a semistar operation on an integral domain $D$. Suppose that
$n:=\widetilde{\star}$-$\dim(D)$. Then
$$n+1\leq\star[X]\text{-}\dim(D[X])\leq 2n+1.$$
\end{thm}

\begin{proof} Consider a chain
$P_1\subseteq\cdots\subseteq P_n$ of quasi-$\widetilde{\star}$-prime
ideals of $D$. Let $Q:=P_n[X]+(X)$. Since $Q\cap D=(P_n[X]+(X))\cap
D=P_n\in\QSpec^{\widetilde{\star}}(D)$, we have, using Remark
\ref{q} that, $Q$ is a quasi-$\star[X]$-prime ideal of $D[X]$. Then
$$P_1[X]\subseteq\cdots\subseteq P_n[X]\subseteq P_n[X]+(X),$$ is a
chain of $n+1$ quasi-$\star[X]$-prime ideals of $D[X]$. Hence
$n+1\leq \star[X]\text{-}\dim(D[X])$.

For the second inequality suppose that $Q\in\QMax^{\star[X]}(D[X])$
is such that $$\hight_{D[X]}Q=\star[X]\text{-}\dim(D[X]).$$ Hence by
\cite[Theorem 38]{K} we obtain that
$\hight_{D[X]}Q\leq2(\hight_D(Q\cap D))+1\leq2n+1$. Consequently we
have $\star[X]\text{-}\dim(D[X])\leq 2n+1$.
\end{proof}

In \cite[Theorem 3]{S2}, Seidenberg showed that for any pair of
positive integers $(n,m)$ with $n+1\leq m\leq2n+1$, there exists a
domain $D$ such that $\dim(D)=d_D$-$\dim(D)=n$ and
$\dim(D[X])=d_{D[X]}$-$\dim(D[X])=d_D[X]$-$\dim(D[X])=m$.

If $X_1,\cdots,X_r$ are indeterminates over $D$, for $r\geq2$, we
let
$$\star[X_1,\cdots,X_r]:=(\star[X_1,\cdots,X_{r-1}])[X_r],$$ where
$\star[X_1,\cdots,X_{r-1}]$ is a stable semistar operation of finite
type on $D[X_1,\cdots,X_{r-1}]$.

\begin{thm}\label{nodim} Let $\star$ be a semistar operation on an integral domain $D$. Suppose that $D$
is a $\widetilde{\star}$-Noetherian domain of
$\widetilde{\star}$-Krull dimension $n$. Then
$$\star[X_1,\cdots,X_m]\text{-}\dim(D[X_1,\cdots,X_m])=n+m.$$
\end{thm}

\begin{proof} Since $D[X_1,\cdots,X_{m-1}]$ is $\star[X_1,\cdots,X_{m-1}]$-Noetherian
domain, it suffices to prove the theorem for the case $m=1$. By
Theorem \ref{dim}, we have $n+1\leq\star[X]$-$\dim(D[X])$. Now let
$M$ be an arbitrary quasi-$\star[X]$-maximal ideal of $D[X]$. Then
$M$ is either an upper to zero, or $P:=M\cap
D\in\QSpec^{\widetilde{\star}}(D)$. Note that in either case $D_P$
is a Noetherian domain (\cite[Proposition 3.8]{EFP}). Hence:
\begin{align*}
\hight_{D[X]}M= & \dim(D[X]_M)=\dim(D_P[X]_{MD_P[X]}) \\[1ex]
              \leq & \dim(D_P[X])=\dim(D_P)+1\\[1ex]
              \leq & n+1.
\end{align*}
The third equality holds since $D_P$ is a Noetherian domain and
\cite[Theorem 30.5]{G}, and the second inequality holds by
\cite[Lemma 2.11]{EFP}. So that by \cite[Lemma 2.11]{EFP} we obtain
that
$$\star[X]\text{-}\dim(D[X])=\sup\{\hight_{D[X]}M|M\in\QMax^{\star[X]}(D[X])\}\leq
n+1,$$ which ends the proof.
\end{proof}

\begin{thm}\label{prudim} Let $\star$ be a semistar operation on an integral domain $D$.
Suppose that $D$ is a P$\star$MD of $\widetilde{\star}$-Krull
dimension $n$. Then $\star[X]\text{-}\dim(D[X])=n+1$.
\end{thm}

\begin{proof} Use the fact that if $D$ is a Pr\"{u}fer
domain, then $\dim(D[X])=\dim(D)+1$ \cite[Corollary]{S2} and by the
same argument as Theorem \ref{nodim} the proof is complete.
\end{proof}

In Corollary \ref{prdim}, we show that if $D$ is a P$\star$MD then
$$\star[X_1,\cdots,X_m]\text{-}\dim(D[X_1,\cdots,X_m])=\widetilde{\star}\text{-}\dim(D)+m.$$

One of the key concepts of Jaffard in \cite{Jaf1}, is that of a
\emph{special chain}, defined as follows. A chain
$\mathcal{C}=\{P_i\}_{i=0}^{m}$ of primes in a polynomial ring
$D[X_1,\cdots,X_m]$ is called a \emph{special chain} if, for each
$P_i\in\mathcal{C}$, the ideal $(P_i\cap D)[X_1,\cdots,X_m]$ is a
member of $\mathcal{C}$. Jaffard's \emph{special chain theorem}
asserts that, if $Q$ is a prime ideal of $D[X_1,\cdots,X_m]$ of
finite height, then $\hight(Q)$ can be realized as the length of a
special chain of primes in $D[X_1,\cdots,X_m]$ with terminal element
$Q$. In particular, if $D$ is a finite dimensional domain, then
$\dim(D[X_1,\cdots,X_m])$ can be realized as the length of a special
chain of prime ideals of $D[X_1,\cdots,X_m]$ (see \cite[Corollary
30.19]{G} for a simple proof). So we make the following remark.

\begin{rem} Let $\star$ be a semistar operation on an integral domain $D$.
If $\widetilde{\star}$-$\dim(D)$ is finite, then
$\star[X_1,\cdots,X_m]\text{-}\dim(D[X_1,\cdots,X_m])$ can be
realized as the length of a special chain of
quasi-$\star[X_1,\cdots,X_m]$-prime ideals of $D[X_1,\cdots,X_m]$.
In fact there exists a quasi-$\star[X_1,\cdots,X_m]$-maximal ideal
$Q$ of $D[X_1,\cdots,X_m]$ such that
$$\star[X_1,\cdots,X_m]\text{-}\dim(D[X_1,\cdots,X_m])=\hight Q.$$
Now by Jaffard's special chain theorem \cite[Corollary 30.19]{G},
$\hight(Q)$ can be realized as the length of a special chain
$(0)=Q_0\subseteq Q_1\subseteq\cdots\subseteq Q_n$ of prime ideals
in $D[X_1,\cdots,X_m]$ with $Q_n=Q$. Since $Q_n$ is a
quasi-$\star[X_1,\cdots,X_m]$-prime ideal of $D[X_1,\cdots,X_m]$,
then each of $Q_1,\cdots, Q_{n-1}$ is a
quasi-$\star[X_1,\cdots,X_m]$-prime ideal of $D[X_1,\cdots,X_m]$ by
Theorem \ref{main}(a) and \cite[Lemma 4.1, and Remark 4.5]{FH}.
\end{rem}

As an application of Theorem \ref{dim} is the following result,
which is the semistar version of \cite[Theorem 8]{S1}.

\begin{thm}\label{dim2} Let $\star$ be a semistar operation on an integral domain $D$.
Suppose that $\widetilde{\star}$-$\dim(D)=1$. Then
$\star[X]$-$\dim(D[X])=2$ if and only if $D$ is a
$\widetilde{\star}$-quasi-Pr\"{u}fer domain.
\end{thm}

\begin{proof} $(\Rightarrow)$. Suppose the contrary. Hence by
\cite[Lemma 2.3]{CF}, there exists an upper to zero $Q$ of $D[X]$
such that $c_D(Q)^{\widetilde{\star}}\varsubsetneq
D^{\widetilde{\star}}$. Then $c_D(Q)^{\widetilde{\star}}$ is
contained in a quasi-$\widetilde{\star}$-prime ideal $P$ of $D$ and
hence $Q\varsubsetneq P[X]$. So that
$2\leq\hight_{D[X]}(P[X])\leq\star[X]\text{-}\dim(D[X])=2$, that is
$\hight_{D[X]}(P[X])=2$. This means that $P[X]$ is a
quasi-$\star[X]$-maximal ideal of $D[X]$. Therefore since
$(P[X]+(X))\cap D=P\in\QSpec^{\widetilde{\star}}(D)$, we obtain that
$P[X]+(X)\in\QSpec^{\star[X]}(D[X])$. Hence $P[X]=P[X]+(X)$ since
$P[X]$ is a quasi-$\star[X]$-maximal ideal of $D[X]$. So that
$(X)\subseteq P[X]$. Consequently $D=c_D((X))\subseteq
c_D(P[X])\subseteq P$, which is a contradiction.

$(\Leftarrow)$. By Theorem \ref{dim} we have
$2\leq\star[X]\text{-}\dim(D[X])\leq3$. If
$\star[X]\text{-}\dim(D[X])=3$, then $\hight_{D[X]}(M)=3$ for some
$M\in\QMax^{\star[X]}(D[X])$. By \cite[Corollary 30.2]{G}, $M$ can
not be an upper to zero. So that $P:=M\cap
D\in\QMax^{\widetilde{\star}}(D)$. From \cite[Lemma 2.1]{CF} and the
hypothesis, we obtain that $D_P$ is a quasi-Pr\"{u}fer domain of
dimension $1$. Hence $\dim(D_P[X])=2$ by \cite[Proposition
30.14]{G}. So we have:
$$3=\hight_{D[X]}(M)=\dim(D[X]_M)=\dim(D_P[X]_{MD_P[X]})\leq\dim(D_P[X])=2,$$
which is a contradiction. Hence $\star[X]\text{-}\dim(D[X])=2$.
\end{proof}

Recall that an integral domain $D$ is called a UM$t$-domain (UM$t$
means ``uppers to zero are maximal $t$-ideals") if every upper to
zero in $D[X]$ is a maximal $t$-ideal \cite[Section 3]{HZ}. It is
observed in \cite[Corollary 2.4 (b)]{CF} that $D$ is a
$w$-quasi-Pr\"{u}fer domain if and only if $D$ is a UM$t$-domain.

\begin{cor} Let $D$ be an integral domain.
Suppose that $w$-$\dim(D)=1$. Then $w[X]$-$\dim(D[X])=2$ if and only
if $D$ is a UM$t$ domain.
\end{cor}

\begin{cor} Let $\star$ be a semistar operation on an integral domain $D$.
Suppose that $\widetilde{\star}$-$\dim(D)=1$. The following
statements are equivalent:
\begin{itemize}
\item[(1)] $D$ is a P$\star$MD.

\item[(2)] $D^{\widetilde{\star}}$ is integrally closed and
$\star[X]$-$\dim(D[X])=2$.
\end{itemize}
\end{cor}

\begin{proof} The equivalence follows easily from Theorem \ref{dim2} and from the fact that
$D$ is a P$\star$MD if and only if, $D$ is a
$\widetilde{\star}$-quasi-Pr\"{u}fer domain and
$D^{\widetilde{\star}}$ is integrally closed, \cite[Lemma 2.17]{CF}.
\end{proof}

In the following result we collect the semistar (Krull) dimension
properties of $[\star]$.

\begin{prop}\label{CF} Let $\star$ be a semistar operation on an integral domain $D$. Suppose that
$n:=\widetilde{\star}\text{-}\dim(D)$. Then
$n\leq[\star]\text{-}\dim(D[X])\leq 2n$. Moreover if $D$ is a
$\widetilde{\star}$-Noetherian domain or a P$\star$MD, then
$[\star]\text{-}\dim(D[X])=\widetilde{\star}\text{-}\dim(D)$.
\end{prop}

\begin{proof} Consider a chain
$P_1\subseteq\cdots\subseteq P_n$ of quasi-$\widetilde{\star}$-prime
ideals of $D$. Since $P_1[X]\subseteq\cdots\subseteq P_n[X]$ is a
chain of $n$ quasi-$[\star]$-prime ideals of $D[X]$, we have $n\leq
[\star]\text{-}\dim(D[X])$. For the second inequality suppose that
$Q\in\QMax^{[\star]}(D[X])$ is such that
$$\hight_{D[X]}Q=[\star]\text{-}\dim(D[X]).$$ If $Q$ is an upper to
zero, then $\hight_{D[X]}Q\leq1\leq2n$. Otherwise by \cite[Theorem
2.3 (e)]{CF1}, there exists a quasi-$\widetilde{\star}$-maximal
ideal $P$ of $D$ such that $Q=P[X]$. Hence by \cite[Theorem 38]{K}
we obtain that $\hight_{D[X]}Q\leq2(\hight_D(P))\leq2n$.
Consequently we have $[\star]\text{-}\dim(D[X])\leq 2n$.

Now suppose that $D$ is a $\widetilde{\star}$-Noetherian domain or a
P$\star$MD. We know that
$\widetilde{\star}\text{-}\dim(D)\leq[\star]$-$\dim(D[X])$. Let $M$
be an arbitrary quasi-$[\star]$-maximal ideal of $D[X]$. Then $M$ is
either an upper to zero, or $M=P[X]$ for some
$P\in\QMax^{\widetilde{\star}}(D)$ by \cite[Theorem 2.3 (e)]{CF1}.
Note that in either case $D_P$ is a Noetherian domain by
\cite[Proposition 3.8]{EFP} (resp. a valuation domain by
\cite[Theorem 3.1]{FJS}). Hence:
\begin{align*}
\hight_{D[X]}P[X]= & \dim(D[X]_{P[X]})=\dim(D_P[X]_{PD_P[X]}) \\[1ex]
              \leq & \dim(D_P[X])-\dim(D_P[X]/PD_P[X])\\[1ex]
              = & \dim(D_P[X])-\dim((D_P/PD_P)[X])\\[1ex]
              = & \dim(D_P)\leq \widetilde{\star}\text{-}\dim(D).
\end{align*}
The fourth equality holds since $D_P$ is a Noetherian domain and
\cite[Theorem 30.5]{G} (resp. a valuation domain and \cite[Theorem
4]{S2}) and the second inequality holds by \cite[Lemma 2.11]{EFP}.
So that by \cite[Lemma 2.11]{EFP} we obtain that
$[\star]\text{-}\dim(D[X])\leq \widetilde{\star}\text{-}\dim(D)$,
which ends the proof.
\end{proof}

Analogous to Seidenberg, in \cite[Theorem 2.10]{F2}, Wang, showed
that for any pair of positive integers $(n,m)$ with $1\leq n\leq
m\leq2n$, there exists a domain $D$ such that $w_D$-$\dim(D)=n$ and
$w_{D[X]}$-$\dim(D[X])=m$. Note that $[w_D]=w_{D[X]}$ by
\cite[Theorem 2.3]{CF1}.

\begin{rem} Let $D$ be an integral domain which is $w_D$-Noetherian
and of $w_D$-dimension $n$. Then
$[w_D]\text{-}\dim(D[X])=w_{D[X]}\text{-}\dim(D[X])=n$ by
Proposition \ref{CF}, while $w_D[X]\text{-}\dim(D[X])=n+1$ by
Theorem \ref{nodim}. This means that $w_D[X]\neq w_{D[X]}(=[w_D])$.
Actually noting Part $(c)$ of Theorem \ref{main}, we have
$w_D[X]\lneq[w_D]$.
\end{rem}

\section{Semistar-valuative dimension}

It is worth reminding the reader of the nice behavior of the
valuative dimension with respect to polynomial rings, in the sense
that $\dim_v(D[X_1,\cdots,X_n])=\dim_v(D)+n$ for each positive
integer $n$ and each ring $D$ (\cite[Theorem 2]{Jaf1}). In this
section we define the \emph{semistar-valuative dimension} of
integral domains and derive its properties.

For this section we need to recall the notion of $\star$-valuation
overring (a notion due essentially to P. Jaffard \cite[page
46]{Jaf}). For a domain $D$ and a semistar operation $\star$ on $D$,
we say that a valuation overring $V$ of $D$ is a
\emph{$\star$-valuation overring of $D$} provided
$F^{\star}\subseteq FV$, for each $F\in f(D)$. Note that, by
definition, the $\star$-valuation overrings coincide with the
$\star_f$-valuation overrings. By \cite[Theorem 3.9]{FL}, $V$ is a
$\widetilde{\star}$-valuation overring of $D$ if and only if $V$ is
a valuation overring of $D_P$ for some quasi-$\star_f$-maximal ideal
$P$ of $D$. Also $V$ is a $\star$-valuation overring of $D$ if and
only if $V^{\star_f}=V$, (cf. \cite[Page 34]{EFP}).

Let $R$ be a B\'{e}zout domain. Then each (nonzero) finitely
generated ideal of $R$ is principal. So that if $J$ is a nonzero
finitely generated ideal of $R$, then $J=J^t$ , and hence each
nonzero ideal of $R$ is a $t$-ideal. This implies that the
$d_R$-operation on $R$ is a unique (semi)star operation of finite
type on $R$. Therefore every (semi)star operation of finite type on
a valuation domain, is the trivial identity operation. The following
result is the key lemma in this section.

\begin{lem}\label{V} Let $\star$ be a semistar operation on an integral domain
$D$. Suppose that $W$ is a valuation overring of $D[X]$. Then $W$ is
a $\star[X]$-valuation overring of $D[X]$, if and only if $W\cap K$
is a $\widetilde{\star}$-valuation overring of $D$.
\end{lem}

\begin{proof} $(\Rightarrow)$. Suppose that $W$ is a $\star[X]$-valuation overring
of $D[X]$. Then by \cite[Theorem 3.9]{FL}, there exists a
$Q\in\QMax^{\star[X]}(D[X])$, such that $D[X]_Q\subseteq W$. Put
$P:=Q\cap D$. Note that $D[X]_Q=D_P[X]_{QD_P[X]}$. Therefore
$D_P[X]\subseteq W$, and hence $D_P\subseteq W\cap K$. If $P=0$,
then $K=W\cap K$, and hence clearly $W\cap K$ is a
$\widetilde{\star}$-valuation overring of $D$. If $P\neq0$, then
$P^{\widetilde{\star}}=(Q\cap D)^{\widetilde{\star}}\varsubsetneq
D^{\widetilde{\star}}$ by Remark \ref{q}. Hence
$P\in\QSpec^{\widetilde{\star}}(D)$. Choose a
quasi-$\widetilde{\star}$-maximal ideal $M$ of $D$ containing $P$ by
\cite[Lemma 2.3 (1)]{FL}. So that $D_M\subseteq D_P\subseteq W\cap
K$. Therefore $W\cap K$ is a $\widetilde{\star}$-valuation overring
of $D$ by \cite[Theorem 3.9]{FL}.

$(\Leftarrow)$. Let $M$ be the maximal ideal of $W$, and set
$Q:=M\cap D[X]$. We need to show that $Q$ is a
quasi-$\star[X]$-prime ideal of $D[X]$. Note that $M\cap K$ is the
maximal ideal of $W\cap K$ by \cite[Theorem 19.16]{G}. Since $W\cap
K$ is a $\widetilde{\star}$-valuation overring of $D$, we have
$(W\cap K)^{\widetilde{\star}}=W\cap K$ by \cite[Page 34]{EFP}. Thus
$\widetilde{\star}_{\iota}$ is a (semi)star operation of finite type
by \cite[Proposition 3.4]{P1}, on $W\cap K$, where $\iota$ is the
canonical inclusion of $D$ to $W\cap K$. So that since $W\cap K$ is
a valuation domain it is the identity operation. Put $P:=Q\cap
D=(M\cap K)\cap D$. If $P=0$ then by construction of $\star[X]$, $Q$
is a quasi-$\star[X]$-prime ideal of $D[X]$. So assume that
$P\neq0$. Now we show that $P^{\widetilde{\star}}\neq
D^{\widetilde{\star}}$. If not
$$
D^{\widetilde{\star}}=P^{\widetilde{\star}}=((M\cap K)\cap
D)^{\widetilde{\star}}=(M\cap K)^{\widetilde{\star}}\cap
D^{\widetilde{\star}}=(M\cap K)\cap D^{\widetilde{\star}}.
$$
Hence $D^{\widetilde{\star}}\subseteq M\cap K$ and therefore,
intersecting with $D$ we find that $D=M\cap D$, which is a
contradiction. Now using Remark \ref{q}, we see that $Q$ is a
quasi-$\star[X]$-prime ideal of $D[X]$. Now choose a
quasi-$\star[X]$-maximal ideal $\mathcal{M}$ of $D[X]$ containing
$Q$. Thus we have $D[X]_{\mathcal{M}}\subseteq D[X]_Q\subseteq W$.
Consequently by \cite[Theorem 3.9]{FL}, we obtain that $W$ is a
$(\widetilde{\star[X]}=)\star[X]$-valuation overring of $D[X]$.
\end{proof}

The following theorem is one of the main results of this section,
whose proof is based on that of \cite[Theorem 30.8]{G}. First, we
need the following definition. Let $D$ be a domain and $T$ an
overring of $D$. Let $\star$ and $\star'$ be semistar operations on
$D$ and $T$, respectively. One says that $T$ is
\emph{$(\star,\star')$-linked to} $D$ (or that $T$ is a
$(\star,\star')${\it -linked overring of} $D$) if
$$F^{\star}=D^{\star}\Rightarrow (FT)^{\star'}=T^{\star'}$$ for each nonzero finitely generated ideal $F$ of
$D$. It was proved in \cite[Theorem 3.8]{EF} that $T$ is
$(\star,\star')$-linked to $D$ if and only if $\Na(D,\star)\subseteq
\Na(T,\star')$.

\begin{thm}\label{Jaffard} Let $\star$ be a semistar operation on an integral domain
$D$, and let $n$ be an integer. Then the following statements are
equivalent:
\begin{itemize}
\item[(1)] Each $(\star,\star')$-linked overring $T$ of $D$ has
$\widetilde{\star'}$-dimension at most $n$, whenever $\star'$ is a
semistar operation on $T$.

\item[(2)] Each $(\star,w_T)$-linked overring $T$ of $D$ has
$w_T$-dimension at most $n$.

\item[(3)] Each overring $T$ of $D$ has
$\widetilde{\star}_{\iota}$-dimension at most $n$, where $\iota:D\to
T$ is the canonical inclusion.

\item[(4)] Each $\widetilde{\star}$-valuation overring of $D$ has
dimension at most $n$.

\item[(5)] For each finite subset $\{t_i\}_{i=1}^n$ of $K$,
$\widetilde{\star}_{\iota}$-$\dim(D[t_1,\cdots,t_n])\leq n$, where
$\iota:D\to D[t_1,\cdots,t_n]$ is the canonical inclusion.

\item[(6)] For each finite subset $\{t_i\}_{i=1}^n$ of $K$,
such that $D[t_1,\cdots,t_n]$ is a $(\star,\star')$-linked overring
of $D$, $\widetilde{\star'}$-$\dim(D[t_1,\cdots,t_n])\leq n$,
whenever $\star'$ is a semistar operation on $D[t_1,\cdots,t_n]$.

\item[(7)]
$\star[X_1,\cdots,X_n]\text{-}\dim(D[X_1,\cdots,X_n])\leq2n$.
\end{itemize}
\end{thm}

\begin{proof} $(1)\Rightarrow(2)$,
$(1)\Rightarrow(3)$, $(1)\Rightarrow(6)$, $(3)\Rightarrow(5)$ and
$(6)\Rightarrow(5)$ are trivial.

$(2)\Rightarrow(4)$. By \cite[Lemma 2.7]{EFP}, $V$ is a ${\widetilde
\star}$-valuation overring of $D$ if and only if $V$ is a
$({\widetilde \star}, d_V)$-linked valuation overring of $D$. The
assertion therefore follows since $w_V=d_V$ for a valuation domain.

$(4)\Rightarrow(3)$. Suppose the contrary. So there exists an
overring $T$ of $D$ containing $P_0\subset P_1\subset\cdots\subset
P_n$, of quasi-$\widetilde{\star}_{\iota}$-prime ideals of $T$,
where $\iota:D\to T$ is the canonical inclusion. Actually one can
choose $P_n$ so that $P_n\in\QMax^{\widetilde{\star}_{\iota}}(T)$.
Consider the chain $P_0T_{P_n}\subset P_1T_{P_n}\subset\cdots\subset
P_nT_{P_n}$ of distinct prime ideals of $T_{P_n}$. Using
\cite[Corollary 19.7]{G}, there exists a valuation overring $V$ of
$T_{P_n}$, such that $V$ contains a chain $M_0\subset
M_1\subset\cdots\subset M_n$ of prime ideals of $V$ and $M_i\cap
T_{P_n}=P_iT_{P_n}$. Since
$P_n\in\QMax^{\widetilde{\star}_{\iota}}(T)$ and $V$ is an overring
of $T_{P_n}$, we obtain that $V$ is a
$\widetilde{\star}_{\iota}$-valuation overring of $T$, by
\cite[Theorem 3.9]{FL}. So that $V^{\widetilde{\star}_{\iota}}=V$,
(see \cite[Page 34]{EFP}). Hence $V^{\widetilde{\star}}=V$.
Therefore $V$ is a $\widetilde{\star}$-valuation overring of $D$
(see \cite[Page 34]{EFP}) and $\dim(V)>n$, which is impossible.

$(5)\Rightarrow(3)$. Suppose there exists an overring $T$ of $D$
containing a chain $P_0\subset P_1\subset\cdots\subset P_n$ of
quasi-$\widetilde{\star}_{\iota}$-prime ideals of $T$, where
$\iota:D\to T$ is the canonical inclusion. Choose $t_i\in
P_i\backslash P_{i-1}$, for each $i=1,\cdots,n$. If
$D'=D[t_1,\cdots,t_n]$, then
$$
(0)\subseteq P_0\cap D'\subset
P_1\cap D'\subset\cdots\subset P_n\cap D'\subset D'.
$$
And since $T$ is an overring of $D'$, $P_0\cap D'\neq0$. Indeed let
$r/s\in P_0$, where $r,s\in D\backslash\{0\}$. Then $r=s(r/s)$ is an
element of $P_0\cap D'$. On the other hand each $P_i\cap D'$ is a
quasi-$\widetilde{\star}_{\iota}$-prime ideals of $D'$, where
$\iota:D\to D'$ is the canonical inclusion. More precisely
$$
(P_i\cap D')^{\widetilde{\star}_{\iota}}\cap D'=(P_i\cap
D')^{\widetilde{\star}}\cap D'=P_i^{\widetilde{\star}}\cap
D'^{\widetilde{\star}}\cap D'=P_i^{\widetilde{\star}}\cap D'=
$$
$$
P_i^{\widetilde{\star}}\cap (T\cap D')=(P_i^{\widetilde{\star}}\cap
T)\cap D'=(P_i^{\widetilde{\star}_{\iota}}\cap T)\cap D'=P_i\cap D'.
$$
Therefore $\widetilde{\star}_{\iota}$-$\dim(D[t_1,\cdots,t_n])> n$,
which is a contradiction.

$(3)\Rightarrow(4)$. Let $V$ be a $\widetilde{\star}$-valuation
overring of $D$. Hence we have $V^{\widetilde{\star}}=V$ by
\cite[Page 34]{EFP}. This means that $\widetilde{\star}_{\iota}$ is
a (semi)star operation on $V$, where $\iota:D\to V$ is the canonical
inclusion. Note that since $\widetilde{\star}_{\iota}$ is of finite
type, then it is the identity operation on the valuation domain $V$.
Thus $\dim(V)=\widetilde{\star}_{\iota}$-$\dim(V)\leq n$.

$(4)\Rightarrow(1)$. Suppose the contrary. So there exists a
$(\star,\star')$-linked overring $T$ of $D$ containing a chain
$P_0\subset P_1\subset\cdots\subset P_n$ of
quasi-$\widetilde{\star'}$-prime ideals of $T$. By the same
reasoning as in the proof of $(4)\Rightarrow(3)$, there exists a
$\widetilde{\star'}$-valuation overring $V$ of $T$ with $\dim(V)>n$.
Thus, by \cite[Lemma 2.7]{EFP}, $V$ is a $(\widetilde{\star'},
d_V)$-linked overring of $T$. Since linked-ness is a transitive
relation (\cite[Theorem 3.8]{EF}), $V$ is a $({\widetilde \star},
d_V)$-linked overring of $D$. Consequently $V$ is a
$\widetilde{\star}$-valuation overring of $D$, which is impossible.

So we showed that $(1)-(6)$ are equivalent.

$(4)\Rightarrow(7)$. To prove
$\star[X_1,\cdots,X_n]\text{-}\dim(D[X_1,\cdots,X_n])\leq2n$, it
suffices in view of what we have just shown, to prove that each
$\star[X_1,\cdots,X_n]$-valuation overring  $W$ of
$D[X_1,\cdots,X_n]$ has dimension at most $2n$. Thus by Lemma
\ref{V}, $W\cap K$ is a $\widetilde{\star}$-valuation overring of
$D$. So that $\dim(W\cap K)\leq n$. Then by \cite[Theorem 20.7]{G},
we have $\dim(W)\leq 2n$.

$(7)\Rightarrow(5)$. We consider a subset $\{t_i\}_{i=1}^n$ of $K$.
If $Q_0$ is the kernel of the $D$-homomorphism
$\varphi:D[X_1,\cdots,X_n]\to D[t_1,\cdots,t_n]$, sending $X_i$ to
$t_i$, then \cite[Lemma 30.7]{G}, shows that $\hight(Q_0)=n$. Note
that $D[t_1,\cdots,t_n]\cong D[X_1,\cdots,X_n]/Q_0$. Suppose that
$\beta\in\QSpec^{\widetilde{\star}_{\iota}}(D[t_1,\cdots,t_n])$ is
such that
$\hight(\beta)=\widetilde{\star}_{\iota}\text{-}\dim(D[t_1,\cdots,t_n]),$
where $\iota:D\to D[t_1,\cdots,t_n]$ is the canonical inclusion.
There exists a prime ideal $Q$ of $D[X_1,\cdots,X_n]$, such that
$\beta=\varphi(Q)\cong Q/Q_0$. We claim that $Q$ is a
quasi-$\star[X_1,\cdots,X_n]$-prime ideal of $D[X_1,\cdots,X_n]$. To
this end set $P:=\beta\cap D$, which, by the same argument as in the
proof of part $(5)\Rightarrow(3)$, is a
quasi-$\widetilde{\star}$-prime ideal of $D$. Note that $Q\cap
D=\beta\cap D=P$. Therefore $(Q\cap D)^{\widetilde{\star}}=
P^{\widetilde{\star}}\subsetneq D^{\widetilde{\star}}$. Then by
repeated applications of Remark \ref{q}, we claim that $Q$ is a
quasi-$\star[X_1,\cdots,X_n]$-prime ideal of $D[X_1,\cdots,X_n]$.
This means that $\hight(Q)\leq2n$ by the hypothesis. Thus we have
$$
\widetilde{\star}_{\iota}\text{-}\dim(D[t_1,\cdots,t_n])=\hight(\beta)=\hight(Q/Q_0)\leq\hight(Q)-\hight(Q_0)\leq2n-n=n,
$$
which ends the proof.
\end{proof}

In \cite{Jaf1} Jaffard defines the \emph{valuative dimension},
denoted $\dim_v(D)$, of the domain $D$ to be the maximal rank of the
valuation overrings of $D$. Now we make the following definition.

\begin{defn} Let $\star$ be a semistar operation on an integral domain
$D$. We say that $D$ has \emph{$\star$-valuative dimension} $n$, and
we write  $\star$-$\dim_v(D)=n$, if each $\star$-valuation overring
of $D$ has dimension at most $n$ and if there exists a
$\star$-valuation overring of $D$ of dimension $n$. If no such
integer exists, we say that the $\star$-valuative dimension of $D$
is infinite.
\end{defn}

Note that $d_D$-$\dim_v(D)=\dim_v(D)$. Since by definition, the
$\star$-valuation overrings coincide with the $\star_f$-valuation
overrings we have $\star_f$-$\dim_v(D)=\star$-$\dim_v(D)$. In
particular $t_D$-$\dim_v(D)=v_D$-$\dim_v(D)$. Suppose that $\star_1$
and $\star_2$ are two semistar operations on an integral domain $D$,
such that $\star_1\leq\star_2$. If $V$ is a $\star_2$-valuation
overring of $D$, then for each $F\in f(D)$ we have
$F^{\star_1}\subseteq F^{\star_2}\subseteq FV$. Hence $V$ is a
$\star_1$-valuation overring of $D$ by definition. So we have:
$$\star_2\text{-}\dim_v(D)\leq\star_1\text{-}\dim_v(D).$$

Using \cite[Corollary 19.7]{G} together with \cite[Theorem 3.9]{FL},
one can easily see that
$\widetilde{\star}$-$\dim(D)\leq\widetilde{\star}$-$\dim_v(D)$. The
following example shows that this inequality is not true in general.

\begin{exam} Let $(D,M)$ be a two dimensional local Noetherian
domain and suppose that $0\subsetneq P\subsetneq M$ be the
corresponding chain of prime ideals. Let $(T_1,N_1)$ and $(T_2,N_2)$
be two rank one discrete valuation rings \cite{Chev} dominating the
local rings $D_P$ and $D$ respectively. Let $\star$ be a semistar
operation on $D$ defined by $E^{\star}=ET_1\cap ET_2$ for each
$E\in\overline{\mathcal{F}}(D)$. Then clearly $\star=\star_f$. We
show that $P, M\in\QSpec^{\star}(D)$. Indeed there exists a positive
integer $k$ such that $PT_1=N_1^k$. Hence $P\subseteq P^{\star}\cap
D=PT_1\cap PT_2\cap D\subseteq PT_1\cap D=N_1^k\cap D\subseteq
N_1\cap D=P$. Therefore $P^{\star}\cap D=P$. By the same way
$M^{\star}\cap D=M$. Therefore we have $\star$-$\dim(D)=2$. Now we
compute $\star$-$\dim_v(D)$. Suppose that $V$ is a $\star$-valuation
overring of $D$. Thus in particular we have $D^{\star}\subseteq DV$
that is $T_1\cap T_2\subseteq V$. Using \cite[Theorem 26.1]{G} we
obtain that $T_1\subseteq V$ or $T_2\subseteq V$. Consequently $\dim
V\leq1$. This means that $\star$-$\dim_v(D)=1$. Thus we have
$$2=\star\text{-}\dim(D)>\star\text{-}\dim_v(D)=1.$$
Note that $\widetilde{\star}=d_D$. So that we have
$\widetilde{\star}\lneq\star$ and
$1=\star\text{-}\dim_v(D)<\widetilde{\star}\text{-}\dim_v(D)=2$.
\end{exam}

By a slight modification of Theorem \ref{Jaffard}, we have:

\begin{thm}\label{J} Let $\star$ be a semistar operation on an integral domain
$D$, and let $n$ be an integer. Then the following statements are
equivalent:
\begin{itemize}
\item[(1)] Each $(\star,\star')$-linked overring $T$ of $D$ has
$\widetilde{\star'}$-dimension at most $n$, and $n$ is minimal,
whenever $\star'$ is a semistar operation on $T$.

\item[(2)] Each $(\star,w_T)$-linked overring $T$ of $D$ has
$w_T$-dimension at most $n$, and $n$ is minimal.

\item[(3)] Each overring $T$ of $D$ has
$\widetilde{\star}_{\iota}$-dimension at most $n$, and $n$ is
minimal, where $\iota:D\to T$ is the canonical inclusion.

\item[(4)] $\widetilde{\star}$-$\dim_v(D)=n$.

\item[(5)] For each finite subset $\{t_i\}_{i=1}^n$ of $K$,
$\widetilde{\star}_{\iota}$-$\dim(D[t_1,\cdots,t_n])\leq n$, and $n$
is minimal, where $\iota:D\to D[t_1,\cdots,t_n]$ is the canonical
inclusion.

\item[(6)] For each finite subset $\{t_i\}_{i=1}^n$ of $K$,
such that $D[t_1,\cdots,t_n]$ is a $(\star,\star')$-linked overring
of $D$, $\widetilde{\star'}$-$\dim(D[t_1,\cdots,t_n])\leq n$, and
$n$ is minimal, whenever $\star'$ is a semistar operation on
$D[t_1,\cdots,t_n]$.

\item[(7)] $\star[X_1,\cdots,X_n]\text{-}\dim(D[X_1,\cdots,X_n])=2n$.
\end{itemize}
\end{thm}

\begin{cor}\label{novdim} Let $\star$ be a semistar operation on an integral domain
$D$. If $D$ is a $\widetilde{\star}$-Notherian domain of
$\widetilde{\star}$-dimension $n$, then
$\widetilde{\star}$-$\dim_v(D)=n$.
\end{cor}

\begin{proof} By Theorem \ref{nodim}, we know
$\star[X_1,\cdots,X_n]\text{-}\dim(D[X_1,\cdots,X_m])=2n$. Hence
$\widetilde{\star}$-$\dim_v(D)=n$.
\end{proof}

Let $D$ be a P$\star$MD. Since for each $M\in\QMax^{\star_f}(D)$,
$D_M$ is a valuation domain by \cite[Theorem 3.1]{FJS}, we have
$\widetilde{\star}$-$\dim(D)=\widetilde{\star}$-$\dim_v(D)$. For an
integer $r$, it is convenient to put $\star[r]$ to denote
$\star[X_1,\cdots,X_r]$ and $D[r]$ to denote $D[X_1,\cdots,X_r]$,
where $X_1,\cdots,X_r$ are indeterminates over $D$.

\begin{cor}\label{asli} Let $\star$ be a semistar operation on an integral domain
$D$. Suppose that $\widetilde{\star}$-$\dim_v(D)=k$. Then
$\star[r]$-$\dim(D[r])=\star[r]$-$\dim_v(D[r])$, for each $r\geq k$.
\end{cor}

\begin{proof} Theorem \ref{J} shows that
$\star[k]$-$\dim(D[k])=2k$. Since $D[r]=D[k][X_{k+1},\cdots,X_r]$,
it follows that:
$$
\star[r]\text{-}\dim(D[r])\geq\dim(D[k])+r-k=2k+r-k=r+k.
$$
If $V$ is a $\star[r]$-valuation overring of $D[r]$, then $V\cap K$
is a $\widetilde{\star}$-valuation overring of $D$ by Lemma \ref{V}.
So that by \cite[Theorem 20.7]{G}, we have $\dim(V)\leq\dim(V\cap
K)+r\leq k+r$. Consequently $\star[r]$-$\dim_v(D[r])\leq k+r$. Since
$\star[r]$-$\dim(D[r])\leq\star[r]$-$\dim_v(D[r])$ is always valid,
we obtain that
$\star[r]$-$\dim(D[r])=\star[r]$-$\dim_v(D[r])=k+r=\widetilde{\star}$-$\dim_v(D)+r$.
\end{proof}

\begin{thm}\label{VV} Let $\star$ be a semistar operation on an integral domain
$D$. Then:
$$\star[m]\text{-}\dim_v(D[m])=\widetilde{\star}\text{-}\dim_v(D)+m.$$
\end{thm}

\begin{proof} Put $n:=\widetilde{\star}\text{-}\dim_v(D)$. If $W$ is a
$\star[m]$-valuation overring of $D[m]$, then by Lemma \ref{V},
$W\cap K$ is a $\widetilde{\star}$-valuation overring of $D$. So
that $\dim(W\cap K)\leq n$. Therefore \cite[Theorem 20.7]{G}, shows
that $\dim(W)\leq n+m$. Consequently
$\star[m]\text{-}\dim_v(D[m])\leq n+m$.

But by assumption, there exists a $\widetilde{\star}$-valuation
overring $V$ of $D$ of rank $n$. So that by \cite[Remark 20.4]{G},
$V$ has an extension to a valuation domain $W$ on
$K(X_1,\cdots,X_m)$, with $\dim(W)=n+m$ and such that
$\{X_1,\cdots,X_m\}$ is contained in the maximal ideal of $W$.
Therefore $W$ is a valuation overring of $D[m]$ of dimension $n+m$.
Since $V=W\cap K$ is a $\widetilde{\star}$-valuation overring of
$D$, Lemma \ref{V} shows that $W$ is a $\star[m]$-valuation overring
of $D[m]$. So that $\star[m]\text{-}\dim_v(D[m])\geq n+m$. Thus we
have
$$\star[m]\text{-}\dim_v(D[m])=n+m,$$
which is the desired equality.
\end{proof}

\begin{cor}\label{sam} Let $\star$ be a semistar operation on an integral domain
$D$. Suppose that
$\widetilde{\star}$-$\dim(D)=\widetilde{\star}$-$\dim_v(D)<\infty$.
If $n$ is a positive integer, then
$$\star[n]\text{-}\dim(D[n])=\star[n]\text{-}\dim_v(D[n])=n+\widetilde{\star}\text{-}\dim(D).$$
\end{cor}

Let $\star$ be a semistar operation on an integral domain $D$.
Recall that the \textit{$\star$-closure} of $D$, defined by:
$$D^{cl^{\star}}:=\bigcup\{(F^{\star}:F^{\star})|F\in f(D)\}$$
is an integrally closed overring of $D$ and, more precisely,
$D^{cl^{\star}}=\bigcap\{V|V\text{ is a }\star\text{-valuation
overring of }D\}$. For more details on this subject and for the
proof of the result recalled above, see \cite{OM1}, \cite{Hal1},
\cite[Proposition 3.2 and Corollary 3.6]{FL1}. Set
$\widetilde{D}:=D^{cl^{\widetilde{\star}}}$ and
$*:=\widetilde{\star}_{\iota}$, where $\iota:D\to \widetilde{D}$ is
the canonical embedding. Note that $\widetilde{*}=*$ by
\cite[Proposition 3.1]{P1}.

\begin{prop}\label{qprdim} Let $\star$ be a semistar operation on an integral domain
$D$. Suppose that $D$ is a $\widetilde{\star}$-quasi-Pr\"{u}fer
domain. Then
$$\widetilde{\star}\text{-}\dim(D)=\widetilde{\star}\text{-}\dim_v(D).$$
\end{prop}

\begin{proof} Recall from \cite[Theorem 2.16]{CF} that, $D$ is a $\widetilde{\star}$-quasi-Pr\"{u}fer
domain if and only if $\Na(\widetilde{D},*)$ is a Pr\"{u}fer domain,
that is $\widetilde{D}$ is a P$*$MD by \cite[Theorem 3.1]{FJS}.
Since $\widetilde{D}$ is a P$*$MD, we have
$*$-$\dim(\widetilde{D})=*$-$\dim_v(\widetilde{D})$. Also an easy
application of \cite[Lemma 2.15]{CF}, yields us that
$\widetilde{\star}$-$\dim(D)=*$-$\dim(\widetilde{D})$. So
$$\widetilde{\star}\text{-}\dim(D)=*\text{-}\dim(\widetilde{D})
=*\text{-}\dim_v(\widetilde{D})=\widetilde{\star}\text{-}\dim_v(D).$$
The last equality holds true since by \cite[Corollary 3.6]{FL1} a
valuation domain is a $\widetilde{\star}$-valuation overring of $D$
if and only if it is a $*$-valuation overring of $\widetilde{D}$
(see Remark \ref{fin} for another reasoning of this equality).
\end{proof}

\begin{cor}\label{prdim} Let $\star$ be a semistar operation on an integral domain
$D$. Suppose that $D$ is a $\widetilde{\star}$-quasi-Pr\"{u}fer
domain (e.g., if $D$ is a P$\star$MD). Then
$$\star[n]\text{-}\dim(D[n])=\star[n]\text{-}\dim_v(D[n])=n+\widetilde{\star}\text{-}\dim(D).$$
\end{cor}

Combining Corollary \ref{prdim} with Theorem \ref{dim2}, we obtain
the following corollary. The special case of $\star=d_D$ is
contained in \cite{S2}.

\begin{cor}\label{dim1=q} Let $\star$ be a semistar operation on an integral domain $D$.
Suppose that $\widetilde{\star}$-$\dim(D)=1$. The following
statements are equivalent:
\begin{itemize}
\item[(1)] $\star[X]$-$\dim(D[X])=2$.

\item[(2)] $\star[m]\text{-}\dim(D[m])=m+1$ for any integer $m$.
\end{itemize}
\end{cor}

In \cite{ABDFK}, to honor Jaffard, the authors defined a domain $D$
to be a \emph{Jaffard domain}, in case $\dim(D)=\dim_v(D)<\infty$.
The class of Jaffard domains contains most of the well-known classes
of finite dimensional rings involved in dimension theory of
commutative rings, such as Noetherian domains, Pr\"{u}fer domains,
universally catenarian domains \cite{BDF}, and stably strong
S-domains \cite{MM, Kab}. As the semistar analogue we define:

\begin{defn} Let $\star$ be a semistar operation on an integral domain
$D$. The domain $D$ is said to be a $\widetilde{\star}$-Jaffard
domain, if $\widetilde{\star}$-$\dim(D)<\infty$ and
$\widetilde{\star}$-$\dim(D)=\widetilde{\star}$-$\dim_v(D)$.
\end{defn}

Note that the notion of $d_D$-Jaffard domain coincides with the
``classical'' notion of Jaffard domain. Note that $D$ is
$\widetilde{\star}$-Jaffard domain if and only if
$\widetilde{\star}$-$\dim(D)<\infty$ and
$\star[r]$-$\dim(D[r])=r+\widetilde{\star}$-$\dim(D)$ for every
$r\in\mathbb{N}$. Indeed let $k=\widetilde{\star}$-$\dim_v(D)$. Then
by Corollaries \ref{VV} and \ref{asli} respectively we have
$$k+\widetilde{\star}\text{-}\dim_v(D)=\star[k]\text{-}\dim_v(D[k])=\star[k]\text{-}\dim(D[k])=k+\widetilde{\star}\text{-}\dim(D).$$
Hence $\widetilde{\star}$-$\dim(D)=\widetilde{\star}$-$\dim_v(D)$.
The converse is true by Corollary \ref{sam}.

Every $\widetilde{\star}$-Noetherian domain and every
$\widetilde{\star}$-quasi-Pr\"{u}fer domain (e.g., every P$\star$MD)
of finite $\widetilde{\star}$-dimension is a
$\widetilde{\star}$-Jaffard domain. As Theorem \ref{dim2} and
Corollary \ref{dim1=q} show, if $\widetilde{\star}$-dimension is
one, then $\widetilde{\star}$-quasi-Pr\"{u}fer domains and
$\widetilde{\star}$-Jaffard domains coincide. For the general case,
we have the following theorem. See also \cite[Theorem 4.3]{Sah} for
several other characterizations of
$\widetilde{\star}$-quasi-Pr\"{u}fer domains. The special case of
$\star=d_D$, of the following theorem is contained in \cite{ACE}.

\begin{thm}\label{qJaf} Let $\star$ be a semistar operation on an integral domain $D$.
Suppose that $\widetilde{\star}$-$\dim(D)$ is finite. Then the
following statements are equivalent:
\begin{itemize}
\item[(1)] $D$ is a $\widetilde{\star}$-quasi-Pr\"{u}fer domain.

\item[(2)] Each $(\star,\star')$-linked overring $T$ of $D$ is a
$\widetilde{\star'}$-quasi-Pr\"{u}fer domain, where $\star'$ is a
semistar operation on $T$.

\item[(3)] Each $(\star,\star')$-linked overring $T$ of $D$ is a
$\widetilde{\star'}$-Jaffard domain, where $\star'$ is a semistar
operation on $T$.

\item[(4)] Each overring $T$ of $D$ is a $\widetilde{\star}_{\iota}$-Jaffard domain, where $\iota$ is
the canonical embedding of $D$ into $T$.
\end{itemize}
\end{thm}

\begin{proof} $(1)\Rightarrow(2)$. Suppose that $D$ is a $\widetilde{\star}$-quasi-Pr\"{u}fer
domain. Hence $\Na(D,\star)$ is a quasi-Pr\"{u}fer domain by
\cite[Theorem 2.16]{CF}. If $T$ is a $(\star,\star')$-linked
overring of $D$, where $\star'$ is a semistar operation on $T$, then
by \cite[Theorem 3.8]{EF}, we have $\Na(D,\star)\subseteq
\Na(T,\star')$. Consequently $\Na(T,\star')$ is a quasi-Pr\"{u}fer
domain by \cite[Corollary 6.5.14]{FHP}. Therefore $T$ is a
$\widetilde{\star'}$-quasi-Pr\"{u}fer domain by \cite[Theorem
2.16]{CF}.

$(2)\Rightarrow(3)$ and $(3)\Rightarrow(4)$ are trivial.

$(4)\Rightarrow(1)$. In order to show that $D$ is a
$\widetilde{\star}$-quasi-Pr\"{u}fer domain, it suffices by
\cite[Theorem 2.16]{CF}, to show that $D_P$ is a quasi-Pr\"{u}fer
domain for all $P\in\QMax^{\widetilde{\star}}(D)$. And for this, it
suffices to prove that each overring $T$ of $D_P$, is a Jaffard
domain by \cite[Theorem 6.7.4]{FHP}. To this end let $P$ be an
arbitrary quasi-$\widetilde{\star}$-maximal ideal of $D$, and $T$ be
an overring of $D_P$. Let $V$ be a valuation overring of $T$. Since
$D_P\subseteq V$, and $P$ is a quasi-$\widetilde{\star}$-maximal
ideal of $D$, we have $V$ is a $\widetilde{\star}$-valuation
overring of $D$ by \cite[Theorem 3.9]{FL}. Thus
$V^{\widetilde{\star}}=V$ by \cite[Page 34]{EFP}. This means that
$V$ is a $\widetilde{\star}_{\iota}$-valuation overring of $T$
(\cite[Page 34]{EFP}), where $\iota$ is the canonical embedding of
$D$ into $T$. So we obtain that
$\dim_v(T)=\widetilde{\star}_{\iota}$-$\dim_v(T)$. Therefore by the
hypothesis we have:
$$\dim(T)\leq\dim_v(T)=\widetilde{\star}_{\iota}\text{-}\dim_v(T)=\widetilde{\star}_{\iota}\text{-}\dim(T)\leq\dim(T).$$
Thus $\dim(T)=\dim_v(T)$, that is $T$ is a Jaffard domain. Hence
$D_P$ is a quasi-Pr\"{u}fer domain for all
$P\in\QMax^{\widetilde{\star}}(D)$, that is $D$ is a
$\widetilde{\star}$-quasi-Pr\"{u}fer domain.
\end{proof}

Recall that if $D$ is a Krull domain then it is a P$v$MD (c.f.
\cite[Remark 4.2]{EFP}). Hence from the above theorem, it can be
seen that a Krull domain is $w$-Jaffard. There is an old question
(see \cite{BK}) asking if is it possible to find a UFD (or a Krull
domain) which is not Jaffard. So, the natural question is the
following: is it possible to find a $w$-Jaffard non Jaffard domain?

Next, we wish to establish that, if $D$ is a
$\widetilde{\star}$-Jaffard domain, then $\Na(D,\star)$ is a Jaffard
domain. First we compute the Krull dimension of the $\star$-Nagata
ring.

\begin{thm}\label{kdim} Let $\star$ be a semistar operation on an integral domain
$D$. Then $\dim(\Na(D,\star))=\star[X]\text{-}\dim(D[X])-1$. In
particular if $D$ is a $\widetilde{\star}$-Jaffard domain, then
$\dim(\Na(D,\star))=\widetilde{\star}\text{-}\dim(D)$.
\end{thm}

\begin{proof} Note that if $Q$ is an upper to zero, then,
$\hight(Q)\leq1$. Also if $Q\in\Spec(D[X])$, and $P:=Q\cap D$, such
that $P[X]\subsetneq Q$, then $\hight(Q)=\hight(P[X])+1$ by
\cite[Lemma 30.17]{G}. So we have:
\begin{align*}
\star[X]\text{-}\dim(D[X])= & \sup\{\hight(Q)|Q\in\QMax^{\star[X]}(D[X])\} \\[1ex]
                      = & \sup\{\hight(Q)|Q\cap D\in\QMax^{\widetilde{\star}}(D)\}\\[1ex]
                      = & \sup\{\hight(P[X])+1|P\in\QMax^{\widetilde{\star}}(D)\}\\[1ex]
                      = & \sup\{\hight(P[X])|P\in\QMax^{\widetilde{\star}}(D)\}+1\\[1ex]
                      = & \dim(\Na(D,\star))+1.
\end{align*}
For the third equality note that if $Q\in\QMax^{\star[X]}(D[X])$,
and $P:=Q\cap D$, then $P[X]\subsetneq Q$. Otherwise $Q=P[X]$. Note
that $P\in\QSpec^{\widetilde{\star}}(D)$ (or equal to zero). Due to
the fact that $(P[X]+(X))\cap D=P$, we obtain by Remark \ref{q} that
$P[X]+(X)\in\QSpec^{\star[X]}(D[X])$. Since
$P[X]\in\QMax^{\star[X]}(D[X])$  and is contained in $P[X]+(X)$, we
have  $P[X]=P[X]+(X)$. Then $(X)\subseteq P[X]$ and therefore
$D=c_D((X))\subseteq c_D(P[X])\subseteq P$ which is a contradiction.
For the last equality note that
$\Max(\Na(D,\star))=\{P\Na(D,\star)|P\in\QMax^{\widetilde{\star}}(D)\}$
\cite[Proposition 3.1 (3)]{FL}.
\end{proof}

Next we compute the valuative dimension of the $\star$-Nagata ring.
Before that, we need some observations and one lemma. Let $D$ be an
integral domain and $\star$ a semistar operation on $D$. One can
consider the contraction map $h:\Spec(\Na(D,\star))\to
\QSpec^{\widetilde{\star}}(D)\cup\{0\}$. Indeed if $N$ is a prime
ideal of $\Na(D,\star)$, then there exists a
quasi-$\widetilde{\star}$-maximal ideal $M$ of $D$, such that
$N\subseteq M\Na(D,\star)$. So that $$h(N)=N\cap D\subseteq
M\Na(D,\star)\cap D=M\Na(D,\star)\cap K\cap
D=M^{\widetilde{\star}}\cap D=M.$$ The third equality holds by
\cite[Proposition 3.4 (3)]{FL}. So that $h(N)\in
\QSpec^{\widetilde{\star}}(D)\cup\{0\}$, since it is contained in
$M$ and \cite[Lemma 4.1 and Remark 4.5]{FH}. Note that if $P\in
\QSpec^{\widetilde{\star}}(D)$, then
$$h(P\Na(D,\star))=P\Na(D,\star)\cap D=P\Na(D,\star)\cap K\cap
D=P^{\widetilde{\star}}\cap D=P.$$ Therefore
$h(\Spec(\Na(D,\star)))=\QSpec^{\widetilde{\star}}(D)\cup\{0\}$. In
fact using \cite[Theorem 2.16]{CF}, the map $h$ is bijective if and
only if $D$ is a $\widetilde{\star}$-quasi-Pr\"{u}fer domain.

\begin{lem}\label{sv} Let $\star$ be a semistar operation on an integral domain
$D$. Then each valuation overring of $\Na(D,\star)$ is a
$\star[X]$-valuation overring of $D[X]$.
\end{lem}

\begin{proof} Let $W$ be a valuation overring of $\Na(D,\star)$. Let $M$
be the maximal ideal of $W$. Set $\fQ:=M\cap\Na(D,\star)$ and
$Q:=M\cap D[X]$. Since $\fQ\in\Spec(\Na(D,\star))$, we have
$h(\fQ)=\fQ\cap D=Q\cap D\in\QSpec^{\widetilde{\star}}(D)\cup\{0\}$.
Thus by Remark \ref{q}, we obtain that $Q$ is a
quasi-$\star[X]$-prime ideal of $D[X]$. Now choose a
quasi-$\star[X]$-maximal ideal $\mathcal{M}$ of $D[X]$ containing
$Q$. Thus we have $D[X]_{\mathcal{M}}\subseteq D[X]_Q\subseteq W$.
Consequently by \cite[Theorem 3.9]{FL}, we obtain that $W$ is a
$(\widetilde{\star[X]}=)\star[X]$-valuation overring of $D[X]$.
\end{proof}

Recall that for each domain $D$,
$\dim_v(D)=\sup\{\dim_v(D_M)|M\in\Max(D)\}$. In fact if
$n=\dim_v(D)$, then there exists a valuation overring $V$, with
maximal ideal $N$, of $D$ such that $\dim(V)=n$. Put $M:=N\cap D$.
So that $V$ is a valuation overring of $D_M$. Hence
$\dim_v(D)=n=\dim_v(V)\leq\dim_v(D_M)\leq\dim_v(D)=n$. Actually one
can assume that $M$ is a maximal ideal of $D$.

Let $\star$ be a semistar operation on an integral domain $D$.
Recall from \cite{FL2} that the \emph{Kronecker function ring of $D$
with respect to the semistar operation $\star$} is defined by:
$$\Kr(D,\star):=\left\{f/g\bigg| \begin{array}{l} f,g\in D[X], g\neq0,\text{ and there exists }h\in D[X]\backslash\{0\}
\\\text{ with } (c(f)c(h))^{\star}\subseteq(c(g)c(h))^{\star} \end{array} \right\}.$$
It is an overring of the $\star$-Nagata ring with quotient field
$K(X)$, which is a B\'{e}zout domain \cite{FL2}. From \cite[Theorem
3.5]{FL1}, we have $V$ is $\star$-valuation overring of $D$ if and
only if $V(X)$ is a valuation overring of $\Kr(D,\star)$. Now we are
ready to prove the following theorem.

\begin{thm}\label{vdim} Let $\star$ be a semistar operation on an integral domain
$D$. Then
$$\widetilde{\star}\text{-}\dim_v(D)=\dim_v(\Na(D,\star)).$$
\end{thm}

\begin{proof} Consider the following inequalities:
\begin{align*}
\widetilde{\star}\text{-}\dim_v(D)\leq\dim_v(\Kr(D,\widetilde{\star}))\leq & \dim_v(\Na(D,\star)) \\[1ex]
                      \leq & \star[X]\text{-}\dim_v(D[X])=\widetilde{\star}\text{-}\dim_v(D)+1.
\end{align*}
The first inequality follows from the fact that if $V$ is a
$\widetilde{\star}$-valuation overring of $D$, then $V(X)$ is a
valuation overring of $\Kr(D,\widetilde{\star})$ and that
$\dim(V)=\dim(V(X))$; second inequality follows from the fact that
$\Na(D,\star)\subseteq \Kr(D,\widetilde{\star})$, while the third
one uses the Lemma \ref{sv}. So that we can assume that
$\widetilde{\star}\text{-}\dim_v(D)$ and $\dim_v(\Na(D,\star))$ are
finite numbers. Now by observation before the theorem, choose a
quasi-$\widetilde{\star}$-maximal ideal $P$ of $D$, such that the
maximal ideal $M:=P\Na(D,\star)$ has the property that
$$\dim_v(\Na(D,\star))=\dim_v(\Na(D,\star)_M)=\dim_v(D_P(X)).$$
But since $P\in\QMax^{\widetilde{\star}}(D)$, each valuation
overring of $D_P$, is a $\widetilde{\star}$-valuation overring of
$D$ \cite[Theorem 3.9]{FL}. Hence we find the inequality
$\dim_v(D_P)\leq\widetilde{\star}\text{-}\dim_v(D)$. Consequently we
have
$$\widetilde{\star}\text{-}\dim_v(D)\leq\dim_v(\Na(D,\star))=\dim_v(D_P(X))=
\dim_v(D_P)\leq\widetilde{\star}\text{-}\dim_v(D),$$ in which the
second equality holds by \cite[Proposition 1.22]{ABDFK}. Thus we
find the desired equality
$\widetilde{\star}\text{-}\dim_v(D)=\dim_v(\Na(D,\star)).$
\end{proof}

As an immediate corollary we have:

\begin{cor} Let $\star$ be a semistar operation on an integral domain
$D$. Then:
\begin{itemize}
\item[(a)] $D[X]$ is a $\star[X]$-Jaffard domain, if and only if, $\Na(D,\star)$ is a
Jaffard domain.

\item[(b)] $D$ is a $\widetilde{\star}$-Jaffard domain if and
only if $\Na(D,\star)$ is a Jaffard domain and
$\star[X]$-$\dim(D[X])=\widetilde{\star}$-$\dim(D)+1$.
\end{itemize}
\end{cor}

\begin{proof} Both statements are easy consequences of Theorems \ref{kdim}
and \ref{vdim}, and for $(a)$ use also Theorem \ref{VV}.
\end{proof}

\begin{rem} By the proof of the above theorem, we have
$\widetilde{\star}\text{-}\dim_v(D)=\dim_v(\Kr(D,\widetilde{\star}))$.
Since $\Kr(D,\widetilde{\star})$ is a B\'{e}zout, and hence a
Pr\"{u}fer domain, we have
$$\widetilde{\star}\text{-}\dim_v(D)=\dim_v(\Kr(D,\widetilde{\star}))=\dim(\Kr(D,\widetilde{\star})).$$
\end{rem}

\begin{rem}\label{fin} Let $D$, $\widetilde{D}$, $\star$, and $*$ be as in the Proposition \ref{qprdim}.
Note that by the proof of part
$(6_{\widetilde{\star}})\Rightarrow(10_{\star_f})$ of \cite[Theorem
2.16]{CF} we have $\Na(\widetilde{D},*)=\overline{\Na(D,\star)}$. So
that by Theorem \ref{vdim} we have

$$*\text{-}\dim_v(\widetilde{D})=\dim_v(\Na(\widetilde{D},*))=\dim_v(\overline{\Na(D,\star)})=\dim_v(\Na(D,\star))=\widetilde{\star}\text{-}\dim_v(D),$$
which is another reason for the last equality in the proof of
Proposition \ref{qprdim}.
\end{rem}

\begin{center} {\bf ACKNOWLEDGMENT}
\end{center}

My sincere thanks goes to Muhammad Zafrullah for his useful comments
on this paper and to Marco Fontana for the valuable suggestions
regarding Remark 3.9 and discussions on the notes after Theorem
4.14.

\end{document}